\documentclass[A4paper,11pt]{article}
\usepackage{latexsym}
\usepackage{mathrsfs}
\usepackage{amssymb}
\usepackage{amscd}
\usepackage[dvips]{graphicx}                  

\newtheorem{theorem}{Theorem}

\newtheorem{lemma}{Lemma}

\newenvironment{proof}[1][Proof]{\textbf{#1.} }{\ \rule{0.5em}{0.5em}}
\long\def\symbolfootnote[#1]#2{\begingroup%
	\def\thefootnote{$\;$}\footnote[#1]{$^*$#2}\endgroup}
\begin{document}
	
	\title{A note on special subsets of the Rudin-Frol\'ik order for regulars}
	\author{Joanna Jureczko\footnote{The author is partially supported by Wroc\l{}aw Univercity of Science and Technology grant of K34W04D03 no. 8211104160.}}
\maketitle

\symbolfootnote[2]{Mathematics Subject Classification: Primary 03E10, 03E20, 03E30..

\hspace{0.2cm}
}

\begin{abstract}
	We show that there is a set of  $2^{2^{\kappa}}$ ultrafilters incomparable in Rudin-Frol\'ik order of $\beta \kappa \setminus \kappa$, where $\kappa$ is regular, for which no subset with more than one element has an infimum.  	
\end{abstract}

\section{Introduction}
Considerations around the Rudin-Frol\'ik order, as shown in the literature, are an important topic, but still little known. This order has been defined by Z. Frol\'ik in 
\cite{ZF} who used it to prove that $\beta\omega\setminus \omega$ is not homogeneous.  M.E. Rudin, who nearly defined this ordering  in \cite{MR1},  as the first observed that the relation between filters she used is really  ordering.
D. Booth in \cite{DB} showed that this relation is a partial ordering of the equivalences classes, that is a tree, and that it is not well-founded.

In \cite{MR} the author defined and studied the partial orders on the type of points in $\beta\omega$ and in $\beta\omega \setminus \omega$. These definitions were used later in \cite{BB} and \cite{BE2}.

E. Butkovi\v cov\'a between 1981 and 1990 published a number of papers concerning ultrafilters  in Rudin-Frol\'ik order in $\beta \omega\setminus \omega$.
In \cite{BB} with L. Bukovsk\'y and in \cite{BE2} she constructed an ultrafilter on $\omega$ with the countable set of its predecessors. In \cite{BE1} she  constructed ultrafilters without immediate predecessors. In \cite{BE3}, Butkovi\v cov\'a showed that there exists in Rudin-Frol\'ik order an unbounded chain  orded-isomorphic in $\omega_1$.  In \cite{BE4}, she proved that  there is a set of $2^{2^{\aleph_0}}$ ultrafilters incomparable in Rudin-Frol\'ik order which is bounded from below and no its subset of cardinality more than one has an infimum. In \cite{BE5}, Butkovi\v cov\'a proved that for every cardinal between $\omega$ and $\mathfrak{c}$ there is a strictly decreasing chain without a lower bound.
In most of these papers there is used method presented in \cite{KK}.

In 1976 A. Kanamori published a paper \cite{AK} in which, among others showed that the Rudin-Frol\'ik tree cannot be very high if one consider it over a measurable cardinal. Moreover, in the same paper he left a number of open problems about Rudin-Frol\'ik order. 
Recently,  M. Gitik in \cite{MG} answered some of them but using metamathematical methods. The solution of some of the problems from \cite{AK} presented in combinatorial methods are in preparation, (\cite{JJ_kanamori}).

However, the Rudin-Frol\'ik order was investigated mainly for $\beta\omega$, significant results may be obtained when considering this order for the space $ \beta \kappa $, where $ \kappa $ is any cardinal. For this purpose there are needed special techniques.

In 2001, Baker and Kunen  presented in \cite{BK} very usefull method which can be recognized as a generalization of method  presented in \cite{KK}. 
It is worth empahsizing that both methods, (from \cite{KK} and \cite{BK}), provide usefull "technology"  for keeping the transfinite construction for an ultrafilter not finished before $\mathfrak{c}$ steps, (see \cite{KK}), and $2^\kappa, $  for $\kappa$ being infinite cardinal, (see \cite{BK}), but the second method has some limitations, among others $\kappa$ must be regular.

As already mentioned, the method from \cite{BK} can be useful in keeping the results for the Rudin-Frol\'ik order but for a regular cardinal $\kappa$.
Due to the lack of adequate useful method for a singular cardinal $\kappa$, the similar results but for singulars are still left as open questions. So far, we have not found an answer whether the assumptions can be omitted, which would probably also involve  changes in the methods used in our considerations. Therefore, based on the results from \cite{BK}, we restrict our results to this particular case.

The results contained in this paper are a continuation of research on the properties of subsets of $\beta\kappa$ (\cite{JJ_order1, JJ_order2}), but due to the methods used here and taken from \cite{BK}, they are also limited to the case where $\kappa$  is a regular cardinal.

The results shown in this paper are suggested by the work of \cite{BE1, BE4} where Butkovi\v cov\'a considered the case of subsets of $ \beta\omega\setminus \omega $ ordered by Rudin-Frol\'ik order. The methods from \cite{BK} allow us to generalize and extend her results to subsets of $ \beta\kappa$, where $ \kappa $ is regular.
The main result of this work is Theorem 1, that there is a set of  $ 2^{2^{\kappa}}$  incomparable ultrafiters in the sense of the Rudin-Frol\'ik order such that although this set is bounded from below, (i.e. there is an ultrafilter which is the predecessor of all ultrafilters of this set), each subset of it of cardinality greater than one has no infimum.

The work is divided into two sections. Section 2 presents the definitions and facts needed later in the paper, but the concept of stratified sets was adopted from \cite{BE1}. The definitions related to the existence of independent matrices are taken from \cite{BK} but appropriately modified to meet current needs.
Section 3 contains the main result along with the auxiliary lemma.

We have tried to present all the necessary definitions, assuming tacitly that the reader has a basic knowledge of ultrafilters and the Rudin-Frol\'ik order.
However, for definitions and facts not quoted here, I refer the reader to e.g. \cite{TJ, CN}.

\section{Preliminaries}

\textbf{2.1.}
In the whole paper, we assume that $\kappa$ is an infinite cardinal. Then $\beta\kappa$ means the \v Cech-Stone compactification, where $\kappa$ has the discrete topology. Hence, $\beta\kappa$ is the space of ultrafilters on $\kappa$ and $\beta \kappa \setminus \kappa$ is the space of nonprincipal ultrafilters on $\kappa$. 
\\\\
\textbf{2.2.}
A set $\{\mathcal{F}_\alpha \colon \alpha < \kappa\}$ of filters on $\kappa$ is \textit{$\kappa$-discrete} iff there is a partition $\{A_\alpha \colon \alpha < \kappa\}$ of $\kappa$ such that $A_\alpha \in \mathcal{F}_\alpha$ for each $\alpha < \kappa$. 
\\\\
\textbf{2.3.} Let $\mathcal{F}, \mathcal{G}\in \beta \kappa\setminus \kappa$. We define \textit{Rudin-Frol\'ik order} as follows
$$\mathcal{F} \leqslant_{RF} \mathcal{G} \textrm{ iff } \mathcal{G} = \Sigma(X, \mathcal{F})$$
for some $\kappa$-discrete set $X = \{\mathcal{F}_\alpha \colon \alpha < \kappa\} \subseteq \beta\kappa,$
where $$\Sigma(X, \mathcal{F}) = \{A \subseteq \kappa \colon \{\alpha < \kappa \colon A \in \mathcal{F}_\alpha\}\in \mathcal{F}\}.$$
We define $$\mathcal{F} =_{RF} \mathcal{G} \textrm{ iff } \mathcal{F} \leqslant_{RF} \mathcal{G} \textrm{ and } \mathcal{G} \leqslant_{RF} \mathcal{F}$$  $$\mathcal{F} <_{RF} \mathcal{G} \textrm{ iff } \mathcal{F} \leqslant_{RF} \mathcal{G} \textrm{ and } \mathcal{F} \not =_{RF} \mathcal{G}.$$
Conversely, if $\mathcal{G} \in \overline{X}$ then there exists a unique ultrafilter $\Omega(X, \mathcal{G})$ such that $\Sigma(X, \Omega(X, \mathcal{G})) = \mathcal{G}$.
\\\\
\textbf{2.4.} The \textit{type} of ultrafilter $\mathcal{F}$ is the set
$$\tau(\mathcal{F}) = \{\mathcal{G} \in \beta \kappa \colon h(\mathcal{F}) = \mathcal{G} \textrm{ for some homeomorphism } h \colon \beta \kappa \to \beta \kappa\}.$$
\textbf{2.5.} 	Let us accept the following notation: 
\begin{itemize}
	\item $\mathcal{FR}(\kappa) = \{A \subset \kappa \colon |\kappa\setminus A|< \kappa\},$
	\item $[A, B, C,...]$ means the filter generated by $A, B, C, ...$.
\end{itemize} 

 An ultrafilter $\mathcal{F}$ in $\beta \kappa$ is called \textit{uniform} if $\mathcal{F} = \{A \subset \kappa \colon |A| = \kappa\}$.
The set of all uniform ultrafilters is denoted by $\mathfrak{u}(\kappa)$.
\\\\
\textbf{Fact 1. [Pospi\v sil]} Let $\kappa$ be a cardinal. Then
$$|\mathfrak{u}(\kappa)| = 2^{2^\kappa}.$$ 
\noindent
\textbf{2.6.} Let $\tau$ and $\kappa$ be infinite cardinals. A set of filters $\{\mathcal{F}_{\xi, \zeta} \colon \xi < \tau, \zeta < \kappa\}$ is \textit{stratified} iff
\begin{itemize}
	\item [(1)] $\{\mathcal{F}_{\xi, \zeta} \colon \zeta < \kappa\}$ is $\kappa$-discrete for each $\xi < \tau$,
	\item [(2)] for each $\xi < \tau, \zeta < \kappa$ and each $\nu$ such that $\xi < \nu < \tau$    $$|\{\mu < \kappa \colon A \in \mathcal{F}_{\nu, \mu}\}| = \kappa$$
	for all $A \in \mathcal{F}_{\xi, \zeta}$.
\end{itemize}
\textbf{2.7.} Let $Y =\{\mathcal{F}_{\xi, \zeta} \colon \xi < \tau, \zeta < \kappa\}$ be a stratified set of filters and let $W$ be a subset of $Y$. We define
\begin{itemize}
	\item [(1)] $W(0) = W$,
	\item [(2)] $W(\gamma) = \bigcup_{\beta < \gamma} W(\beta)$ for limit $\gamma$,
	\item [(3)] $W(\gamma+1) = W(\gamma) \cup \{\mathcal{F}_{\xi, \zeta} \colon \exists_{\eta> \gamma}\  \exists_{{A \in \mathcal{F}_{\xi, \zeta}}} \{\mathcal{F}_{\eta, \nu}\colon A \in \mathcal{F}_{\eta, \nu}\}\subseteq W(\gamma)\}$,
	\item [(4)] $\tilde{W} = \bigcup_{\gamma < \kappa^+} W(\gamma)$.
\end{itemize}

Intuitively, the above construction is used to select only certain filters from $Y$ with the desired property, (see e.g. property (P) in 2.9), and then add (inductively) to the set $ W $ only those filters outside $ W $ which satisfy the condition (2) in the definition in 2.6. This construction   will be used to define property (P), (see 2.9), the formulation of which would not be possible taking the entire set of ultrafiters into account.
\\\\
\textbf{2.8.} A function $\hat\varphi \colon [\kappa^+]^{<\omega} \to [\kappa]^{< \omega}$ is \textit{$\kappa$-shrinking} iff 
\begin{itemize}
	\item [(1)] $p\subseteq q$ implies $\hat{\varphi}(p) \subseteq \hat{\varphi}(q)$, for any $p, q \in [\kappa^+]^{<\omega}$,
	\item [(2)] $\hat{\varphi}(0) = 0$.
\end{itemize}

A \textit{step family} (over $\kappa$, with respect to $\hat{\varphi}$) is a family of subsets of $\kappa$, $$\{E_t \colon t \in [\kappa]^{<\omega}\} \cup \{A_\alpha \colon \alpha < \kappa^+\}$$ satisfying the following conditions:
\begin{itemize}
	\item [(1)] $E_s\cap E_t = \emptyset$ for all $s, t \in [\kappa]^{<\omega}$ with $s \not =t$,
	\item [(2)] $|\bigcap_{\alpha \in p} A_\alpha \cap  \bigcup_{t\not \supseteq \hat{\varphi}(p)}E_t| < \kappa$ for each $p \in [\kappa^+]^{<\omega}$,
	\item [(3)] if $\hat{\varphi}(p) \subseteq t$, then $|\bigcap_{\alpha\in p} A_\alpha \cap E_t|=\kappa$ for each $p \in [\kappa^+]^{<\omega}$ and $t \in [\kappa]^{<\omega}$. 
\end{itemize}

Let $I$ be an index set and $\mathcal{F}$ be a filter on $\kappa$. The family $$\{E_t^i \colon t \in [\kappa]^{<\omega}, i \in I\} \cup \{A_\alpha^i \colon \alpha < \kappa^+, i \in I\}$$ is an \textit{independent matrix of $|I|$ step-families} (over $\kappa$) with respect to $\mathcal{F}, \hat{\varphi}$ iff
\begin{itemize}
	\item [(1)] for each fixed $i \in I$, $\{E_t^i \colon t \in [\kappa]^{<\omega}\} \cup \{A_\alpha^i \colon \alpha < \kappa^+\}$ is a step-family,
	\item [(2)] if $n \in \omega, p_0, p_1, ..., p_{n-1} \in [\kappa^+]^{<\omega}, t_0, t_1, ..., t_{n-1} \in [\kappa]^{<\omega}$,  $i_0, i_1, ..., i_{n-1} \in I$ with $i_k\not = i_m, k\not = m$ and $\hat{\varphi}(p_k) \subseteq t_k$, then 
	$$\kappa \setminus \bigcap_{k=1}^{n-1}(\bigcap_{\alpha\in p_k}A^{i_k}_{\alpha} \cap E^{i_k}_{t_k}) \not \in \mathcal{F}.$$  
\end{itemize}
\noindent
\textbf{Fact 2 (\cite{BK}).} If $\kappa$ is a regular cardinal and $\hat{\varphi}$ is a $\kappa$-shrinking function, then there exists an independent matrix of $2^\kappa$ step-families over $\kappa$ with respect to the  filter $\mathcal{FR}(\kappa)$, $\hat{\varphi}$. 
\\\\
\textbf{2.9.} Let $\{D_t \colon t < [\kappa]^{<\omega}\}$ be a partition of $\kappa$ and let $\hat{\varphi}$ be a $\kappa$-shrinking function. A stratified set of filters $\{\mathcal{F}_{\xi, \zeta} \colon \xi < \tau, \zeta < \kappa\}$ \textit{satisfies property (P)} iff $\mathcal{F}_{\nu, \mu} \not \in \tilde{W}$
implies that there exists  $\{W_\eta \colon \eta < \kappa^+\} \subseteq \mathcal{F}_{\nu, \mu}$ such that for all $p \in [\kappa^+]^{<\omega}$
$$ |\bigcap_{\eta \in p} W_\eta \cap D_{\hat{\varphi}(p)}| < \kappa,$$
where $$W = \{\mathcal{F}_{\xi, \zeta} \colon \exists_{t < [\kappa]^{<\omega}} D_t \in \mathcal{F}_{\xi, \zeta}\}.$$
\textbf{2.10.} 
A set of utrafilters $\{\mathcal{F}_{\xi, \zeta} \colon \xi < \tau, \zeta < \kappa\}$ is called \textit{a stratified set with uniform predecessor} $\mathcal{F}$ iff
\begin{itemize}
	\item [(1)] $\{\mathcal{F}_{\xi,\zeta} \colon \xi < \tau, \zeta < \kappa\}$ is stratified,
	\item [(2)] $\{\mathcal{F}_{\xi,\zeta} \colon \xi < \tau, \zeta < \kappa\}$ has the property $(P)$ for each partition of $\kappa$ and each $\kappa$-shrinking function $\{\mathcal{F}_{\xi,\zeta} \colon \xi < \tau, \zeta < \kappa\}$,
	\item [(3)] $\Omega(X_{\xi+1}, \mathcal{F}_{\xi, \zeta}) = \mathcal{F}$ for each $\xi < \tau, \zeta < \kappa$,  where
	$X_\xi = \{\mathcal{F}_{\xi, \zeta} \colon \zeta< \kappa\}$,
	\item [(4)] $\Omega(\bigcup_{\mu< \xi}X_\mu, \mathcal{F}_{\xi, \zeta}) = \mathcal{F}$ for limit $\xi$.
\end{itemize}

\section{Main results}

We start with the simple lemma whose proof we leave to the reader.

\begin{lemma}
	Let $\mathcal{F}$ be an ultrafilter and let $\{D_t \colon t \in [\kappa]^{<\omega}\}$  be a partition of $\kappa$. Let $\mathcal{F}_{t}, t \in [\kappa]^{< \omega}$ denote any ultrafilter extending the filter
	$$ [\mathcal{FR}(\kappa), \{D_t\}]$$ and Let $\mathcal{G}$ denote any ultrafilter extending the filter 
	$$[\mathcal{FR}(\kappa), \{\bigcup_{t \in [\kappa]^{<\omega}} D_t \colon D_t \in \mathcal{F}\} ].$$
	If $\mathcal{G} \in \overline{\{\mathcal{F}_t \colon t \in [\kappa]^{<\omega}\}}$, then $$\Omega(\{\mathcal{F}_t \colon t \in [\kappa]^{<\omega}\}, \mathcal{G}) = \mathcal{F},$$ (i.e. $\Sigma (\{\mathcal{F}_t \colon t \in [\kappa]^{<\omega}\}, \mathcal{F}) = \mathcal{G}$). 
\end{lemma}

The main result (Theorem 1) is  based on the next lemma.

\begin{lemma}
	Let $\kappa$ be a regular cardinal and let $\hat{\varphi}$ be a $\kappa$-shrinking function and let $\mathcal{F}$ be a minimal ultrafilter. Then
	\begin{itemize}
		\item [(i)] There exists a stratified set of ultrafilters $\{\mathcal{F}_{\xi, \zeta} \colon \xi, \zeta < \kappa\}$ with uniform predecessor $\mathcal{F}$.
		\item [(ii)] There is a family $\{\{\mathcal{F}^h_{\xi, \zeta} \colon \xi, \zeta < \kappa\} \colon h \in 2^{2^\kappa}\}$ of stratified sets of ultrafilters with uniform predecessor $\mathcal{F}$ such that if $h \not = h', (h, h'\in 2^{2^\kappa})$, then there exists $C \subseteq \kappa$ such that for each $\xi, \zeta < \kappa$ $$C \in \mathcal{F}^h_{\xi, \zeta} \textrm{ and } \kappa \setminus C \in \mathcal{F}^{h'}_{\xi, \zeta}.$$	
	\end{itemize}
	\end{lemma} 

\begin{proof}
	We start with proving $(i)$. By Fact 2, fix a matrix $$\{E_t^i \colon t \in [\kappa]^{<\omega}, i \in 2^\kappa\} \cup \{A_\eta^i \colon \eta < \kappa^+, i \in 2^\kappa\}$$ of $2^\kappa$ step-families over $\kappa$ independent with respect to the filter $\mathcal{FR}(\kappa), \hat \varphi$.
	For our purpose, we slightly modify this matrix by "shrinking" $A^i_\alpha$ to $A^i_\alpha \subseteq \bigcup \{E_t^i \colon t \in [\kappa]^{<\omega}\}$ and then expanding $E^i_t$ so that $\{E^i_t \colon t \in [\kappa]^{< \omega}\}$ is a partition of $\kappa$. 
	Thus, we obtain the matrix fullfilling the following conditions:
	\begin{itemize}
		\item [(a)] $E^i_s\cap E^i_t = \emptyset$ for all $s, t \in [\kappa]^{<\omega}$ with $s \not =t$,
		\item [(b)] $|\bigcap_{\eta \in p} A^i_\eta \cap  \bigcup_{t\not \supseteq \hat{\varphi}(p)}E^i_t| < \kappa$ for each $p \in [\kappa^+]^{<\omega}$,
		\item [(c)] if $\hat{\varphi}(p) \subseteq t$, then $|\bigcap_{\eta\in p} A^i_\eta \cap E^i_t|=\kappa$ for each $p \in [\kappa^+]^{<\omega}$ and $t \in [\kappa]^{<\omega},$ 
		\item [(d)] $\bigcup\{E^i_t \colon t \in [\kappa]^{< \omega}\} = \kappa,$
		\item [(e)] $|\bigcap_{\eta \in p} A^i_\eta \setminus  \bigcup_{t \supseteq \hat{\varphi}(p)}E^i_t| < \kappa$ for each $p \in [\kappa^+]^{<\omega}.$
	\end{itemize}
Note that $(b)$ and $(d)$ implies $(e)$.
Moreover, the condition $(b)$ is still preserved after expanding $\{E_t^i \colon t \in [\kappa]^{<\omega}\}$ to a partition of $\kappa$. 
\\Indeed. If there are $p_0 \in [\kappa^+]^{<\omega}$ and $i_0 \in 2^\kappa$ such that
$$|\bigcap_{\eta \in p_0} A_\eta^{i_0} \cap  \bigcup_{t\not \supseteq \hat{\varphi}(p_0)}E_t^{i_0}| = \kappa,$$
then $|\bigcup_{t\not \supseteq \hat{\varphi}(p_0)}E_t^{i_0}| = \kappa.$
Then, by $(d)$ and $(a)$, there would exist  $t_0 \supseteq \hat{\varphi}(p_0)$ such that  $|E_{t_0}^{i_0}| < \kappa$. Hence
$$|\bigcap_{\eta\in p_0} A^{i_0}_\eta \cap E^{i_0}_{t_0}|<\kappa.$$
which contradicts $(c)$. 
	
	Let $\{Z_\alpha \colon \alpha < 2^\kappa\}$ be a family of all subsets of $\kappa$.
	Let $\mathcal{A} = \{\mathcal{A}^\alpha \colon \alpha < 2^\kappa\}$
	be a sequence enumerating all partitions of $\kappa$ in such a way that each partition occurs $2^\kappa$ many times.
	
	Now, we will construct a stratified set $\{\mathcal{F}_{\xi, \zeta} \colon \xi, \zeta < \kappa\}$ of ultrafilters with  uniform predecessor $\mathcal{F}$.
In order to do this we will define $\{B_{\xi, \zeta} \colon \xi, \zeta < \kappa\}$ such that $B_{\xi, \zeta} \cap B_{\xi, \zeta'} = \emptyset$ for $\zeta \not = \zeta'$.

To do this fix a partition
$G = \{G_\zeta \colon \zeta < \kappa, |G_\zeta|=\kappa\}$ of $\kappa$.

For any $\zeta < \kappa$ we will define $B_{0,\zeta}$ as follows.
Let $p \in [\kappa^+]^{<\omega}$ be such that $\zeta \in p$ and denote such $p$ by $p_\zeta$.
Fix $t \in [\kappa]^{<\omega}$ such that $t \supseteq \hat{\varphi}(p_\zeta)$. Such $t$ denote by $t_\zeta$.
Take
$$B_{0,\zeta} = (\bigcap_{\eta \in p_\zeta}A^0_\eta \cap E^0_{t_\zeta}) \setminus (\bigcup_{\mu < \zeta}\bigcap_{\eta \in p_\mu }A^0_\eta\cap E^0_{t_\mu}).$$
By the above construction $B_{0,\zeta}, \zeta < \kappa$ are pairwise disjoint and  well defined, (which follows from $(a)-(e)$).

Assume that we have constructed $B_{\delta, \zeta}$ for some $\delta< \xi< \kappa$.
\\
If $\xi$ is limit, then set $B_{\xi, \zeta} = \bigcap_{\delta < \xi} B_{\delta, \zeta}$.
\\
Now we will construct $B_{\xi, \zeta}$ for $\xi = \delta+1$.
Observe that since $G$ is a partition of $\kappa$ then $\zeta \in G_\gamma$ for some $\gamma < \kappa$. Then take
$$B_{\xi,\zeta} =  B_{\delta, \gamma}\cap [(\bigcap_{\eta \in p_\zeta}A^\xi_\eta \cap E^\xi_{t_\zeta}) \setminus (\bigcup_{\mu < \zeta}\bigcap_{\eta \in p_\mu }A^\xi_\eta\cap E^\xi_{t_\mu})].$$
Thus, for any $\xi < \kappa$ we have constructed the family $\{B_{\xi, \zeta} \colon \zeta < \kappa\}$ of the required property.

Fix a minimal ultrafilter $\mathcal{F}$.
Take a minimal ultrafilter $\mathcal{G}_\zeta$ of the same type as $\mathcal{F}$ such that $G_\zeta \in \mathcal{G}_\zeta$.
For any $\xi, \zeta < \kappa$ define
 $$\mathcal{F}^0_{\xi, \zeta} =[\mathcal{FR}(\kappa), \{B_{\xi, \zeta}\}, \{\bigcup_{\delta \in A} B_{\nu, \delta} \colon A \in \mathcal{G}_\zeta, \nu > \xi\}].$$
 	Let $I_0 = 2^\kappa \setminus \kappa$.
 It is easy observation that  
 $$\{\mathcal{F}^0_{\xi, \zeta} \colon \xi, \zeta < \kappa\}$$
 is a stratified set of filters.
and
 $$\{E_t^i \colon t \in [\kappa]^{<\omega}, i \in I_0\} \cup \{A_\alpha^i \colon \alpha < \kappa^+, i \in I_0\}$$
 is an independent matrix of $|I_0|$ step-families (over $\kappa$) with respect to $\mathcal{F}^0_{\xi,\zeta}, \hat{\varphi}$, for all $\xi, \zeta < \kappa$. 
		
	Now, we will construct filters $\mathcal{F}^\alpha_{\xi, \zeta}$ and indexed sets $I_\alpha$ by induction on $\alpha < 2^\kappa$ steps fulfilling the properties
\begin{itemize}
	\item [(1)] $\mathcal{F}^0_{\xi, \zeta}$ and $I_0$ as are done above,
	\item [(2)] $\mathcal{F}^\alpha_{\xi, \zeta}$ is a filter on $\kappa$,  $I_\alpha\subset 2^\kappa$ and the matrix $$\{E_t^i \colon t \in [\kappa]^{<\omega}, i \in I_\alpha\} \cup \{A_\alpha^i \colon \alpha < \kappa^+, i \in I_\alpha\}$$ of remaining step-families is  independent w.r.t $\mathcal{F}^\alpha_{\xi, \zeta}, \hat{\varphi}$,
	\item [(3)] $\mathcal{F}^\alpha_{\xi, \zeta} = \bigcup_{\beta < \alpha} \mathcal{F}^\beta_{\xi,\zeta}$, $I_\alpha = \bigcap_{\beta<\alpha} I_\beta$, for limit $\alpha$, 
	\item [(4)] $\mathcal{F}^\beta_{\xi, \zeta} \subseteq \mathcal{F}^\alpha_{\xi, \zeta}$, $I_\beta \supseteq I_\alpha$, whenever $\beta < \alpha$,
	
	\item [(5)] $I_\alpha \setminus I_{\alpha + 1}$ is finite,
	\item [(6)] if $\alpha \equiv 0\ (mod\  2)$, then either $Z_\alpha \in \mathcal{F}^\alpha_{\xi,\zeta}$ or $\kappa \setminus Z_\alpha \in \mathcal{F}^\alpha_{\xi, \zeta}$,
	\item [(7)] if $\alpha \equiv 1\ (mod\  2)$,   then the set  $\{\mathcal{F}^\alpha_{\xi, \zeta}\colon  \xi, \zeta < \kappa\}$ of filters is stratified with uniform predecessor $\mathcal{G}_\zeta$.
\end{itemize}
Then, take  $$\mathcal{F}_{\xi, \zeta} = \bigcup_{\alpha < 2^\kappa} \mathcal{F}^\alpha_{\xi, \zeta}.$$ 
Thus, $\{\mathcal{F}_{\xi, \zeta} \colon \xi, \zeta < \kappa\}$ 
will be the required set of ultrafilters.

Assume that $\mathcal{F}^\beta_{\xi, \zeta}$ and $I_\beta$ have been constructed for some $\beta <\alpha < 2^\kappa$. The limit step is done.
We show how to obtain $\mathcal{F}^{\alpha}_{\xi, \zeta}$ and $I_{\alpha}$ for $\alpha = \beta+1$.
\\
\\
\underline{Case $\alpha\equiv 0\ (mod \ 2)$}.

If
$\mathcal{R} = [\mathcal{F}^\beta_{\xi, \zeta}, Z_\beta]$ is a proper filter and the matrix of step-families
$$\{E_t^i \colon t \in [\kappa]^{<\omega}, i \in I_\beta\} \cup \{A_\eta^i \colon \eta < \kappa^+, i \in I_\beta\}$$ is independent w.r.t. $\mathcal{R}, \hat{\varphi}$, then  put $\mathcal{F}^\alpha_{\xi, \zeta}= [\mathcal{F}^\beta_{\xi, \zeta}, Z_\beta]$ and $I_\alpha = I_\beta$.

Otherwise, fix $n \in \omega$, distinct $i_k \in I_\alpha$ and $\hat{\varphi}(p_k)\subseteq t_k,$ for $ k < n$, such that
$$\kappa \setminus [Z_\alpha \cap \bigcap_{k=0}^{n-1}(\bigcap_{\eta\in p_k}A^{i_k}_{\eta} \cap E^{i_k}_{t_k})]\in \mathcal{F}^\alpha_{\xi, \zeta}.$$
Then, put $$\mathcal{F}^{\alpha}_{\xi, \zeta} = [\mathcal{F}^{\beta}_{\xi, \zeta}, \{A^{i_k}_{p_k} \colon 0 \leqslant k \leqslant n-1\}, \{E^{i_k}_{t_k} \colon 0 \leqslant k \leqslant n-1\}]$$ $$I_{\alpha} = I_{\beta}\setminus \{i_k \colon 0 \leqslant k \leqslant n-1\}.$$
Then $\kappa \setminus Z_\alpha \in \mathcal{F}^{\alpha}_{\xi, \zeta}.$

To show that $(2)$ is fulfilled for $\mathcal{F}^{\alpha}_{\xi, \zeta}$ and $I_\alpha$ it is enough to observe that each element of $\mathcal{F}^{\alpha}_{\xi, \zeta}$ is of the form 
$$A \cap (\bigcap_{\eta \in p_k}A^{i_k}_\eta \cap E^{i_k}_{t_k})$$
for some $A \in \mathcal{F}^{\alpha}_{\xi, \zeta} $ and $ 0 \leqslant k \leqslant n-1$. Then,
$$\bigcap_{\eta \in p_k}A^{i_k}_\eta \cap E^{i_k}_{t_k} \supseteq\bigcap_{k=0}^{n-1}(\bigcap_{\eta\in p_k}A^{i_k}_{\eta} \cap E^{i_k}_{t_k}) $$
and $$\bigcap_{k=0}^{n-1}(\bigcap_{\eta\in p_k}A^{i_k}_{\eta} \cap E^{i_k}_{t_k}) \in \mathcal{F}^{\beta}_{\xi, \zeta}.$$
Thus, $(2)$ for $\mathcal{F}^{\alpha}_{\xi, \zeta}$ and $I_\alpha$ is fulfilled which follows from $(2)$ for $\mathcal{F}^{\beta}_{\xi, \zeta}$ and $I_\beta$. 
\\
\\
\underline{Case $\alpha\equiv 1\ (mod \ 2)$}.

If $D^\beta_s \in \mathcal{F}^\beta_{\xi, \zeta}$ for some $s \in [\kappa]^{<\omega}$, then put
$\mathcal{F}^\alpha_{\xi, \zeta} = \mathcal{F}^\beta_{\xi, \zeta}$ and $Z_\alpha = I_\beta$.

Otherwise, we have (by (6))
$\kappa\setminus D^\beta_s \in \mathcal{F}^\beta_{\xi, \zeta}$ for each $s \in [\kappa]^{<\omega}$.

Consider sequences 
$\langle V^\beta_s \colon s \in [\kappa]^{<\omega}\rangle$
such that $V^\beta_t \subseteq V^\beta_s$ whenever $t \subseteq s$ and 
$D^\alpha_s \subseteq V^\alpha_s$ for any $s \in [\kappa]^{<\omega}$.

Now, choose $i \in I_\beta$ such that 
$$\bigcup\{E^i_t\setminus V^i_s \colon s, t \in [\kappa]^{<\omega}, t \subseteq s\} \not = \emptyset.$$
Set
$$W^\beta_\delta = A^i_\delta \cap \bigcup\{E^i_t\setminus V^i_s \colon s, t \in [\kappa]^{<\omega}, t \subseteq s\}.$$
Then, set
$I_\alpha = I_\beta \setminus\{i\}$ and $\mathcal{F}^\alpha_{\xi, \zeta} = [\mathcal{F}^\beta_{\xi, \zeta}, \{W^\beta_\delta \colon \delta < \kappa^+\}].$

To show that $(2)$ holds for $\mathcal{F}^\alpha_{\xi, \zeta}$ and $I_\alpha$ it is enough to observe that each element of $\mathcal{F}^\alpha_{\xi, \zeta}$  is of the form 
$$A \cap \bigcap_{\delta \in p}W^\alpha_\delta$$
for some $A\in \mathcal{F}^\alpha_{\xi, \zeta}$ and $p \in [\kappa^+]^{<\omega}$.
But 
$$\bigcap_{\delta\in p} W^\alpha_\delta = \bigcap_{\delta\in p} A^i_\delta \cap \bigcup\{E^i_t\setminus V^i_s \colon s, t \in [\kappa]^{<\omega}, t \subseteq s\} \cap (\kappa\setminus D^\alpha_{\hat{\varphi}(p)})$$
and $\kappa\setminus D^\alpha_{\hat{\varphi}(p)} \in \mathcal{F}^\beta_{\xi, \zeta}$.
Thus, $(2)$ for $\mathcal{F}^{\alpha}_{\xi, \zeta}$ and $I_\alpha$ is fulfilled which follows from $(2)$ for $\mathcal{F}^{\beta}_{\xi, \zeta}$ and $I_\beta$. 

Using the similar argument, it is easy to check that  $\{\mathcal{F}^\alpha_{\xi, \zeta} \colon \xi, \zeta < \kappa\}$ is a  stratified set of filters.

Now, we show that $\{\mathcal{F}^\alpha_{\xi, \zeta} \colon \xi, \zeta < \kappa\}$ fulfills property $(P)$.
Indeed, for each $p \in [\kappa^+]^{<\omega}$ we have
$$\bigcap_{\delta \in p} W^\alpha_\delta = \bigcap_{\delta\in p} A^i_\delta \cap \bigcup\{E^i_t\setminus V^i_s \colon s, t \in [\kappa]^{<\omega}, t \subseteq s\} $$ 
(by $(e)$)
$$\subseteq^* \bigcup_{t \supseteq\hat{\varphi}(p)}E^i_t \cap \bigcup\{E^i_t\setminus V^i_s \colon s, t \in [\kappa]^{<\omega}, t \subseteq s\}  $$
(by $(a)$)
$$= \bigcup_{t \supseteq \hat{\varphi}(p)} E^i_t \setminus V^\alpha_t \subseteq \bigcup_{t \supseteq \hat{\varphi}(p)} \kappa\setminus V^\alpha_t $$
by monotonicity of $\langle V^\beta_s \colon s \in [\kappa]^{<\omega}\rangle$
$$\subseteq \kappa \setminus V^\alpha_{\hat{\varphi}(p)} \subseteq \kappa \setminus D^\alpha_{\hat{\varphi}(p)}.$$
Thus $|\bigcap_{\delta\in p}W^\alpha_\delta \cap D^\alpha_{\hat{\varphi}(p)}| < \kappa$.
The proof of $(i)$ is complete.
\\

To show that $(ii)$ holds,  it is enough to use the method presented above for different minimal ultrafilters and apply Fact 1. 
\end{proof}
\\

Now, we are ready to prove the main result.

\begin{theorem}
	Let $\kappa$ be a regular cardinal and let $\hat{\varphi}$ be a $\kappa$-shrinking function. Then there exists a set $A \subseteq \beta \kappa \setminus \kappa$ of ultrafilters in Rudin-Frol\'ik order such that
	\begin{itemize}
		\item [(i)] $|A| = 2^{2^\kappa}$
		\item[(ii)] $\exists_{\mathcal{F}} \forall_{\mathcal{G} \in A} \mathcal{F} < \mathcal{G}$
		\item [(iii)] $\forall_{S \subset A} |S|>1$ $\inf S$ does not exist.
	\end{itemize}
	\end{theorem}

\begin{proof}
	By Lemma 1, Lemma 2 and Fact 1  we have the existence of $2^{2^{\kappa}}$ distinct stratified sets 
	$$\{\{\mathcal{F}^h_{\xi, \zeta} \colon \xi, \zeta < \kappa\} \colon h \in 2^{2^\kappa}\}$$
	with uniform predecessor $\mathcal{F}$ being minimal in Rudin-Frol\'ik order.
	
To complete the proof we need to prove two claims.
	\\\\
	\textbf{Claim 1} Let $\mathcal{G}$ be an ultrafilter such that
	such that $\tau(\mathcal{G}) \not = \tau(\mathcal{F})$ and
		$$\mathcal{G} <_{RF} \mathcal{F}^h_{\xi, \zeta}$$   for any $h \in 2^{2^\kappa}$.
Then, there exists an ultrafilter $\mathcal{K}$ which is  an immediate predecessor of $\mathcal{G}$.
\\

\begin{proof} (of Claim 1)
	Let $\mathcal{G}$ be as in Claim. Observe that there exists a $\kappa$-discrete set $W \subseteq \{\mathcal{F}^h_{\xi, \zeta} \colon \xi, \zeta < \kappa\}$ such that $$\mathcal{F}^h_{\xi, \zeta} = \Sigma (W, \mathcal{G}).$$
	By property $(P)$ we have that $\mathcal{F}^h_{\xi, \zeta} \in \tilde{W}$, but  since $\tau(\mathcal{G}) \not = \tau(\mathcal{F})$ we have      $\mathcal{F}^h_{\xi, \zeta} \not \in W(1)$. Observe that $\mathcal{F}^h_{\xi, \zeta} \in \overline{W(1)\setminus W}$. 
	\\Indeed. 
	We proceed by induction.
	
	Assume that for $\beta < \gamma$ we have 
	$$\mathcal{F}^h_{\xi, \zeta} \in W(\beta) \setminus W(1) \textrm{ implies } \mathcal{F}^h_{\xi, \zeta} \in \overline{W(1) \setminus W}.$$
	
	The case when for limit cardinal is obvious.	We show case $\gamma  = \beta + 1$.
	Let $\mathcal{F}^h_{\xi, \zeta} \in W(\gamma) \setminus W(1)$. Then, by induction step $\mathcal{F}^h_{\alpha, \mu} \in \overline{W(1) \setminus W}$ for each $\mathcal{F}^h_{\eta, \mu} \in W(\beta) \setminus W(1)$, $\mu < \kappa$ and some $\alpha > \eta$.
	Then also $\mathcal{F}^h_{\xi, \zeta} \in \overline{W(1)\setminus W}$.

	Thus, $$\mathcal{F}^h_{\xi, \zeta} \in \overline{(W(1)\setminus W)} \setminus W(1)= \overline{W(1)\cap (\tilde{W}\setminus W)})\setminus W(1).$$
	Hence
	$$\Omega(W(1), \mathcal{F}^h_{\xi, \zeta}) <_{RF} \mathcal{G}.$$
	Hence, there exists a $\kappa$-discrete set $$Y=\{\mathcal{F}_\gamma \colon \tau(\mathcal{F}_\gamma) =  \tau( \mathcal{F}), \gamma < \kappa\}$$
	of minimal ultrafilters
	such that 
	$\mathcal{G} = \Sigma(Y, \Omega(W(1), \mathcal{F}^h_{\xi, \zeta}))$.
	
Hence $\mathcal{K} = \Omega(W(1), \mathcal{F}^h_{\xi, \zeta})$ is the required ultrafilter.
\end{proof}
\\\\
\textbf{Claim 2}
Let $\mathcal{G}$ be an ultrafilter such that 
 $$\mathcal{G} <_{RF} \mathcal{F}^h_{\xi, \zeta} \in \{\mathcal{F}^h_{\xi, \zeta} \colon \xi, \zeta < \kappa\}.$$
Then,  there exists an ultrafilter $\mathcal{K}$ such that 
 $$\mathcal{G} <_{RF} \mathcal{K}<_{RF} \mathcal{F}^h_{\xi, \zeta}$$
 for each $h \in 2^{2^\kappa}$.
 \\
 
 \begin{proof}(of Claim 2)
 	Let $\mathcal{G}$ be as in Claim.
 	Then, there is $Y_h \subseteq \{\mathcal{F}^h_{\xi, \zeta} \colon \xi, \zeta < \kappa\}$ such that 
 	$$\mathcal{F}^h_{\xi, \zeta} = \Sigma(Y_h, \mathcal{G}).$$
 	Let $\{D^h_t \colon t \in [\kappa]^{<\omega}\}$ be a partition of $\kappa$ such that 
 	$D^h_t \in \mathcal{F}^h_{\xi, \zeta}$, where $\mathcal{F}^h_{\xi, \zeta} \in Y_h$.
 	Let $$Z_h = \bigcup_{t \in [\kappa]^{<\omega}}\{\mathcal{F}^h_{\nu, \mu}\colon  D^h_t \in \mathcal{F}^h_{\nu, \mu}, \mu > \zeta\}.$$
 	Then $$\mathcal{G} = \Omega(Y_h, \mathcal{F}^h_{\xi, \zeta})< \Omega(Z_h, \mathcal{F}^h_{\xi, \zeta}).$$
 	Observe that 
 	$\Omega(Z_h, \mathcal{F}^h_{\xi, \zeta}) $ is of the same type as $\Sigma(S, \mathcal{G})$, where $S$ is the set of all ultrafilters of the same type as $\mathcal{F}$.
 	Then $\mathcal{K} = \Sigma(S, \mathcal{G})$ is the common predecessor of $\mathcal{F}^h_{\xi, \zeta}$ and is greater than $\mathcal{G}$.
 	
 	By the same argument as above no subset of $A$ of cardinality at least 2 has an infimum. 
 \end{proof}
\\
\\	
	\textbf{The end  of the  proof of Theorem 1.} Now,  it is enough to observe that by Claim 1 each predecessor of $\mathcal{F}^h_{\xi, \zeta}$ has an immediate predecessor, by Claim 2 there exists $A \subseteq \beta \kappa \setminus \kappa$ composed of $2^{2^\kappa}$ incomparable ultrafilters with the common predecessor $\mathcal{F}$ and without a greatest common predecessor.
	Thus, if $S \subseteq A$ and $|S|>1$, then $\inf S$ does not exist.
	\end{proof}
\\\\
\textbf{Acknowledgments.} The author is very grateful to the anonymous reviewer for his insight in reading thie previous version of this paper. Their remarks undoubtedly avoided many inaccuracies and made the text more readable.

\noindent
{\sc Joanna Jureczko}
\\
Wroc\l{}aw University of Science and Technology, Wroc\l{}aw, Poland
\\
{\sl e-mail: joanna.jureczko@pwr.edu.pl}

\end{document}